\documentclass[12pt,a4paper]{amsart}
\usepackage{amsmath}
 \usepackage{amsaddr}
\usepackage{mathabx}
\usepackage{amstext}
\usepackage{amssymb,mathtools}
\usepackage{array}
\usepackage{amsthm}
\usepackage{dsfont}
\usepackage[latin1]{inputenc}
\usepackage[T1]{fontenc}
\usepackage{mathtools}
\usepackage{appendix}
\usepackage{xcolor}
\usepackage{ulem}
\usepackage[arrow, matrix, curve]{xy}
\usepackage{verbatim}
\usepackage{graphicx}
\usepackage{subfigure}
\usepackage{hyperref}
\usepackage{cases}
\usepackage{esvect}
\usepackage[update,prepend]{epstopdf}

\newtheorem{theorem}{Theorem}
\theoremstyle{plain}
\newtheorem{definition}[theorem]{Definition}
\newtheorem{assumption}[theorem]{Annahme}
\newtheorem{lemma}{Lemma}

\newtheorem{proposition}[theorem]{Proposition}
\theoremstyle{definition}
\newtheorem{example}{Example}

\numberwithin{equation}{section}

\def\eps{\varepsilon}

\def\ri{{\rm i}}

\def\bi{\begin{itemize}}
\def\ei{\end{itemize}}

\newcommand{\R}{\mathbb{R}}
\newcommand{\C}{\mathbb{C}}
\newcommand{\B}{\mathbb{B}}

\newcommand{\N}{\mathbb{N}}

\newcommand{\cF}{{\mathcal F}}

\renewcommand{\phi}{\varphi}

\newcommand{\pa}{\partial}

\newcommand{\Bvechat}{\widehat{\vec{B}}}
\newcommand{\Nvechat}{\widehat{\vec{N}}}
\newcommand{\dvechat}{\widehat{\vec{d}}}
 \def\dd{\, {\rm d}}

\DeclareMathOperator{\sign}{sign}
\DeclareMathOperator{\Real}{Re}
\DeclareMathOperator{\Imag}{Im}

\DeclareOldFontCommand{\it}{\normalfont\itshape}{\mathit}

\newcommand{\bspm}{\left(\begin{smallmatrix}}\newcommand{\espm}{\end{smallmatrix}\right)}
\newcommand{\bpm}{\begin{pmatrix}}\newcommand{\epm}{\end{pmatrix}}

\def\bs{\begin{satz}}\def\es{\end{satz}}
\def\blem{\begin{lemma}}\def\elem{\end{lemma}}
\def\bthm{\begin{theorem}}\def\ethm{\end{theorem}}
\def\bcor{\begin{corollary}}\def\ecor{\end{corollary}}
\def\beq{\begin{equation}}\def\eeq{\end{equation}}
\def\beqq{\begin{equation*}}\def\eeqq{\end{equation*}}
\def\bal{\begin{align}}\def\eal{\end{align}}
\def\bpf{\begin{proof}}\def\epf{\end{proof}}
\def\bex{\begin{example}}\def\eex{\end{example}}
\def\brem{\begin{remark}}\def\erem{\end{remark}}
\def\bass{\begin{assumption}}\def\eass{\end{assumption}}
\def\bprop{\begin{proposition}}\def\eprop{\end{proposition}}
\def\bdefi{\begin{definition}}\def\edefi{\end{definition}}

\DeclareFontEncoding{FMS}{}{}
\DeclareFontSubstitution{FMS}{futm}{m}{n}
\DeclareFontEncoding{FMX}{}{}
\DeclareFontSubstitution{FMX}{futm}{m}{n}
\DeclareSymbolFont{fouriersymbols}{FMS}{futm}{m}{n}
\DeclareSymbolFont{fourierlargesymbols}{FMX}{futm}{m}{n}
\DeclareMathDelimiter{\VERT}{\mathord}{fouriersymbols}{152}{fourierlargesymbols}{147}

\def\bi{\begin{itemize}}
\def\ei{\end{itemize}}
\def\ben{\begin{enumerate}}
\def\een{\end{enumerate}}
\usepackage{enumitem}
\setitemize{leftmargin=*}

\begin{document}

\title[Bifurcation of  Gap Solitons in CME\MakeLowercase{s} in $\R^\MakeLowercase{d}$]{Bifurcation of  Gap Solitons in Coupled Mode Equations in $d$ Dimensions}
\author{Tom\'a\v{s} Dohnal$^{1}$ and Lisa Wahlers$^{2}$}

\address{$^{1}$ Martin-Luther-Universit\"at Halle-Wittenberg, Institut f\"ur Mathematik, D-06099 Halle (Saale), Germany\\ 
$^{2}$ Technische Universit\"at Dortmund, Fakult\"at f\"ur Mathematik, Vogelpothsweg 87, D-44227 Dortmund, Germany}
\email{tomas.dohnal@mathematik.uni-halle.de}
\date{\today}


\begin{abstract}
We consider a system of first order coupled mode equations in $\R^d$ describing the envelopes of wavepackets in nonlinear periodic media. Under the assumptions of a spectral gap and a generic assumption on the dispersion relation at the spectral edge, we prove the bifurcation of standing gap solitons of the coupled mode equations from the zero solution. The proof is based on a Lyapunov-Schmidt decomposition in Fourier variables and a nested Banach fixed point argument. The reduced bifurcation equation is a perturbed stationary nonlinear Schr\"odinger equation. The existence of solitary waves follows in a symmetric subspace thanks to a spectral stability result. A numerical example of gap solitons in $\R^2$ is provided.
\end{abstract}

\maketitle

   \vspace*{2mm} {\bf Key-words:} coupled mode equations, gap soliton, bifurcation, fixed point

     \vspace*{2mm}
		
\setcounter{page}{1}

\section{Introduction}
\label{S:intro}

First order coupled mode equations (CMEs) are used to describe a class of wavepackets in periodic structures \cite{AW89,SU01,GWH01,AP05,GMS08,DH17,DW20}. They are modulation equations for the envelopes of asymptotically broad and small wavepackets. They were first studied in nonlinear optical fiber gratings, see e.g. \cite{AW89,AP05}. A rigorous justification of such an approximation was performed in \cite{SU01} and \cite{GWH01} for the one dimensional cubic nonlinear wave equation. In \cite{DW20} the authors derived and justified CMEs as modulation equations for the $d$-dimensional periodic Gross-Pitaevskii equation. These CMEs have the form
\beq
\begin{aligned}
\ri (\pa_tA_j+v_g^{(j)}\cdot \nabla A_j) + \sum_{r=1}^N \kappa_{jr} A_r + N_j(\vec{A})=0, \ j = 1,\dots,N,
\label{E:CME}
\end{aligned}
\eeq
where for $j,r\in \{1,\dots,N\}$
\begin{align*}
N_j(\vec{A})&:=\sum_{(m,n,o)\in\{1,\dots,N\}^3} \gamma_j^{(m,n,o)}A_{m}\overline{A}_{n}A_{o},\\
\gamma^{(m,n,o)}_j & \in \C, \kappa_{jr}\in \C, v_g^{(j)}\in \R^d, 
\end{align*}
and where the matrix $\kappa=(\kappa_{jr})_{j,r=1}^N$ is Hermitian. This system (although only for the setting with $\kappa=0$) was first derived in \cite{GMS08}. Like with all modulation equations, the application of CMEs is not limited to the Gross-Pitaevskii equation. It applies to wavepackets centered around $N$ Bloch waves ($N$-wave mixing) with nonzero group velocities in models with nonlinearities that are cubic at lowest order.

The aim of this paper is to prove the existence of localized time harmonic solutions 
\beq\label{E:GS}
\vec{A}(x,t)=e^{-\ri \omega t} \vec{B}(x), \ |\vec{B}(x)|\to 0 \text{ as } |x|\to \infty,
\eeq
with $\omega$ in a gap of the linear spatial operator of \eqref{E:CME}. Such solutions are often called (standing) gap solitons. A necessary condition for the existence of a spectral gap is $\kappa\neq 0$. Hence, gap solitons cannot be obtained in the setting of \cite{GMS08}.
We prove the existence of gap solitons in an asymptotic region near a spectral edge point $\omega_0$. The result can be interpreted as a bifurcation from the zero solution at the spectral edge.

The equation for $\vec{B}$ is
\beq\label{E:B-eq}
\omega \vec{B}- L(\nabla)\vec{B}+\vec{N}(\vec{B})=0, 
\eeq
where 
$$L(\nabla)=\begin{pmatrix} 
-\ri v_g^{(1)}\cdot \nabla & & \\
& \ddots & \\
& &  -\ri v_g^{(N)}\cdot \nabla 
\end{pmatrix} - \kappa.$$

Our proof is constructive in that we use an asymptotic approximation of a solution $\vec{B}$ at $\omega = \omega_0+O(\eps^2), \eps \to 0$. The approximation is a modulation ansatz with a slowly varying envelope modulating the bounded linear solution at the spectral edge $\omega_0$. The envelope is shown to satisfy a $d-$dimensional nonlinear Schr\"odinger equation (NLS) with constant coefficients. We prove that for sufficiently smooth \textit{PT} symmetric (parity time symmetric) solutions of the NLS there are solutions $\vec{B}$ of \eqref{E:B-eq} at $\omega = \omega_0+O(\eps^2)$ which are close to the asymptotic ansatz. The proof is carried out in Fourier variables in $L^1(\R^d)$. It is based on a decomposition of the solution in Fourier variables according to the eigenvectors of $L(\ri k)\in \C^{n\times n}$ and on a nested Banach fixed point argument. The reduced bifurcation equation is a perturbed stationary nonlinear Schr\"odinger equation (NLS). Solitary waves are then found via a persistence argument starting from solitary waves of the unperturbed NLS. The persistence holds in a symmetric subspace thanks to a spectral stability result of Kato.

The chosen approach is similar to that used in \cite{DPS09,DU09,DP20}. Unlike in these papers, where $L^2$-based spaces were used, we work here in $L^1$ in Fourier variables. This avoids the unfavorable scaling property of the $L^2$ norm of functions with an asymptotically slow dependence on $x$, namely $\|f(\eps\cdot)\|_{L^2(\R^d)}=\eps^{-d/2}\|f\|_{L^2(\R^d)}$. The $L^1-$approach was first used in this context for the bifurcation of time harmonic gap solitons in the one dimensional wave equation in \cite{PS07}.

The question of the existence of solitary waves of CMEs has previously been addressed only in one dimension in \cite{AW89}, where an explicit family of gap solitons was found for CMEs describing the asymptotics of wavepackets in media with infinitesimally small contrast. These gap solitons are parametrized by the velocity $v\in (-1,1)$ (after a rescaling). In \cite{D14} a numerical continuation was used to construct gap solitons also in one dimensional CMEs for finite contrast periodic structures.  It was shown in \cite{DW20}  that for \eqref{E:CME}  in dimensions $d>1$ a spectral gap of $L(\nabla)$ does not exist for $N$ ungerade and for $N=2$. Next, a gap was found in a special case of \eqref{E:CME} with $d=2, N=4$ and standing gap solitons were computed numerically for this case. Here we assume the presence of a spectral gap and prove the existence of standing gap solitons of the form \eqref{E:GS} for $\omega$ asymptotically close to the spectrum under the condition that the spectral edge is given by an isolated extremum of the dispersion relation.

The rest of the paper consists firstly of a formal derivation of the effective NLS equation for the modulation ansatz in Section \ref{S:asymp}. Next, in Section \ref{S:main-thm} we state and prove the main approximation result. Finally, Section \ref{S:num} presents a numerical example of a solution $\vec{B}$ and a numerical verification of the convergence of the asymptotic error.

\section{Formal Asymptotics of Gap Solitons}\label{S:asymp}
The formal asymptotics of localized solutions of \eqref{E:B-eq} were performed already in \cite{DW20}. We repeat here the calculation for readers' convenience.

The spectrum of $L(\nabla)$ can be determined using Fourier variables. We employ the Fourier transform
$$\widehat{f}(k):=(\mathcal{F} f)(k):=(2\pi)^{-d/2}\int_{\R^d}f(x)e^{-\ri k\cdot x}\dd x$$
with the inverse formula $f(x)=(\mathcal{F}^{-1} \widehat{f})(x)=(2\pi)^{-d/2}\int_{\R^d}\widehat{f}(k)e^{\ri k\cdot x}\dd k$. The spectrum of $L(\nabla)$ is 
$$\sigma(L(\nabla))=\cup_{j \in\{1,\dots,N\}}\lambda_j(\R^d),$$ 
where $\lambda_j(k)$ is the eigenvalue of $L(\ri k)\in \C^{N\times N}$ for each $k\in \R^d$, i.e.
$$L(\ri k) \vec{\eta}^{(j)}(k)= \lambda_j(k)\vec{\eta}^{(j)}(k)$$
for some $\vec{\eta}^{(j)}(k) \in \C^N\setminus\{0\}$. Because $L(\ri k)^*=L(\ri k)$ for all $k\in \R^d$, we have $\lambda_j:\R^d\to \R$. The mapping $k\mapsto (\lambda_1(k),\dots,\lambda_N(k))^T$ is the dispersion relation of \eqref{E:CME}.

The central assumptions of our analysis are 
\bi
\item[(A.1)] The spectrum $\sigma (L(\nabla))\subset \R$ has a gap, denoted by $(\alpha,\beta)$ with $\alpha<\beta$.
\item[(A.2)] $\omega_0\in\{\alpha,\beta\}$ and for some $j_0 \in \N, k_0 \in \B$ we have
$$\omega_0=\lambda_{j}(k) \quad \text{if and only if } (j,k)=(j_0,k_0).$$
\ei
As mentioned in the introduction, assumption (A.1) implies that if $d\geq 2$, then $N\geq 4$ and $N$ even. Assumption (A.2) means that the spectral edge $\omega_0$ is defined by one isolated extremum of the eigenvalue $\lambda_{j_0}$ and that this is separated at $k=k_0$ from all other eigenvalues.

We make the following asymptotic ansatz for a gap soliton at $\omega = \omega_0 +\eps^2\omega_1\notin \sigma(L(\nabla))$, where $\eps>0$ is a small parameter and $\omega_1=O(1)$ (as $\eps \to 0$),
\beq\label{E:ans-x}
\vec{B}_{\text{app}}(x):=\eps C(\eps x)e^{\ri k_0\cdot x}\vec{\eta}^{(j_0)}(k_0).
\eeq
In Fourier variables this is 
\beq\label{E:ans}
\Bvechat_{\text{app}}(k)=\eps^{1-d}\widehat{C}\left(\frac{k-k_0}{\eps}\right)\vec{\eta}^{(j_0)}(k_0).
\eeq
Substituting \eqref{E:ans} and $\omega=\omega_0+\eps^2\omega_1$ into the Fourier transform of the left hand side of \eqref{E:B-eq}, we get 
$$
\begin{aligned}
&(\omega_0+\eps^2\omega_1)\Bvechat_{\text{app}} - L(\ri k)\Bvechat_{\text{app}}+\Nvechat(\vec{B}_\text{app})\\
=&\eps(\omega_0+\eps^2\omega_1-\lambda_{j_0}(k))\eta^{(j_0)}(k_0) \widehat{C}\left(\frac{k-k_0}{\eps}\right)+\eps^3 \vec{N}(\eta^{(j_0)}(k_0))(\widehat{C}*\widehat{\overline{C}}*\widehat{C})\left(\frac{k-k_0}{\eps}\right)\\
=&\eps^3\left[\left(\omega_1 - \left(\frac{k-k_0}{\eps}\right)^TG_0\frac{k-k_0}{\eps}\right)\widehat{C}\left(\frac{k-k_0}{\eps}\right) + \Gamma(\widehat{C}*\widehat{\overline{C}}*\widehat{C})\left(\frac{k-k_0}{\eps}\right)\right]\vec{\eta}^{(j_0)}(k_0) - \vec{R}(k),
\end{aligned}
$$
where
$$G_0:=\frac{1}{2}D^2 \lambda_{j_0}(k_0), \quad \Gamma:=\vec{\eta}^{(j_0)}(k_0)^*\Nvechat(\vec{\eta}^{(j_0)}(k_0)),$$
(with $v^*$ being the Hermitian transpose of a vector $v\in \C^N$) and where
$$
\begin{aligned}
\vec{R}(k):=&\eps^3 \left(\lambda_{j_0}(k)-\omega_0-\left(\frac{k-k_0}{\eps}\right)^TG_0\frac{k-k_0}{\eps}\right)\vec{\eta}^{(j_0)}(k_0)\widehat{C}\left(\frac{k-k_0}{\eps}\right) \\
+& \eps^3 \left( \Nvechat(\vec{\eta}^{(j_0)}(k_0)) - \vec{\eta}^{(j_0)}(k_0)^*\Nvechat(\vec{\eta}^{(j_0)}(k_0))\vec{\eta}^{(j_0)}(k_0))\right)(\widehat{C}*\widehat{\overline{C}}*\widehat{C})\left(\frac{k-k_0}{\eps}\right)
\end{aligned}
$$
is small as shown in Sec. \ref{S:main-thm}.

A necessary condition for the smallness of the residual corresponding to $\vec{B}_{\text{app}}$ is the vanishing of the square brackets. This is equivalent to
\beq\label{E:NLS}
\omega_1 C + \nabla^T(G_0 \nabla C) + \Gamma |C|^2C=0
\eeq
for $C:\R^d\to\C$. Equation \eqref{E:NLS} is the effective nonlinear Schr\"odinger equation (NLS) for the envelope $C$.


\section{The Bifurcation and Approximation Result}
\label{S:main-thm}

Under assumptions (A.1-A.2) and the following assumption (A.3) we prove the bifurcation result below.
\bi
\item[(A.3)] The kernel of the Jacobian $J$ corresponding to the NLS equation, as defined in \eqref{E:J}, is $(n+1)$-dimensional (i.e. generated only by the continuous invariances of the NLS).
\ei
We define next the space $L^1_s(\R^d)$ for $s\geq 0$ as
$$L^1_s(\R^d):=\{f\in L^1(\R^d): (1+|\cdot|)^s f \in L^1(\R^d))\}.$$ 
For vector valued functions $f:\R^d\to \C^N$ we write $f \in L_s^1(\R^d)$ if $f_j\in L_s^1(\R^d)$ for each $j=1,\dots,N$.

The space of continuous functions $f:\R^d\to \C$ satisfying the asymptotics $f(x)\to 0$ as $|x|\to\infty$ is denoted by $C_0(\R^d)$. We equip the space with the supremum norm.

\bthm\label{T:main}
Choose $\omega_0$ such that (A.1) and (A.2) are satisfied. Let $(\alpha,\beta)\subset \R$ be the spectral gap from (A.1) and let  $\omega_1\in \R$ be such that $\text{sign}(\omega_1)=1$ if $\omega_0=\alpha$ and $\text{sign}(\omega_1)=-1$ if $\omega_0=\beta$. If $C$ is a $PT$-symmetric (i.e. $C(-x)=\overline{C(x)}$) solution of \eqref{E:NLS}
with $\widehat{C}\in L^1_{4}(\R^d)$ and such that (A.3) holds, then there are constants $c_1,c_2,\eps_0>0$ such that for each $\eps \in (0,\eps_0)$ there is a solution $\vec{B}$ of equation \eqref{E:B-eq} with $\omega=\omega_0+\eps^2\omega_1$ which satisfies $\widehat{\vec{B}}\in L^1_{2}(\R^d)$ and 
$$\left\|\widehat{\vec{B}}-\eps^{1-d}\widehat{C}\left(\frac{\cdot-k_0}{\eps}\right)\vec{\eta}^{(j_0)}(k_0)\right\|_{L^1(\R^d)}\leq c_1\eps^2.$$
In particular,
$$\|\vec{B}-\eps C(\eps \cdot)\vec{\eta}^{(j_0)}(k_0)e^{\ri k_0\cdot}\|_{C_0(\R^d)}\leq c_2\eps^2.$$
The constants $c_1$ and $c_2$ depend polynomially on $\|\widehat{C}\|_{L^1_4(\R^d)}$.
\ethm
Clearly, due to $\widehat{\vec{B}}\in L^1(\R^d)$ the lemma of Riemann-Lebesgue implies the decay $\vec{B}(x)\to 0$ as $|x|\to \infty$.

The existence of a $PT$-symmetric solution $C$ is satisfied, e.g., if $G_0$ is definite and $\text{sign}(\Gamma)=-\sign(\omega_1)$. Due to the extremum of $\lambda_{j_0}$ at $k=k_0$ we have then that $\Gamma$ is positive/negative if $D^2\lambda{j_0}(k_0)$ is positive/negative definite  respectively. Hence, the NLS is of focusing type and after a rescaling of the $x$ variables it supports a real, positive, radially symmetric solution with exponential decay at infinity (Townes soliton).

Note that the condition on $\text{sign}(\omega_1)$ implies
$$\omega=\omega_0+\eps^2 \omega_1 \notin \sigma (L(\nabla)).$$

We proceed with the proof of Theorem \ref{T:main}. Like in Sec. \ref{S:asymp} we work here in Fourier variables. We employ a Lyapunov-Schmidt-like decomposition. For each $k\in\R^d$ we split the solution $\widehat{\vec{B}}(k)\in \C^N$ into the component proportional to the eigenvector $\vec{\eta}^{(j_0)}(k)$ and the $l^2(\C^N)-$orthogonal complement. We define the projections
$$P_k:\C^N\to \text{span}{\vec{\eta}^{(j_0)}(k)}, \vec{v} \mapsto (\vec{\eta}^{(j_0)}(k)^*\vec{v})\vec{\eta}^{(j_0)}(k)$$
and
$$Q_k:=I-P_k.$$
Then 
$$\Bvechat(k) = \Bvechat_P(k)+\Bvechat_Q(k),$$
where
$$\Bvechat_P(k)=P_k\Bvechat(k)=:\psi(k)\vec{\eta}^{(j_0)}(k), \ \Bvechat_Q(k)=Q_k\Bvechat(k).$$
We aim to construct a solution $\Bvechat$  with $\psi$ approximated by the envelope in our ansatz, i.e. by $\eps^{1-d}\widehat{C}\left(\tfrac{\cdot-k_0}{\eps}\right)$. We choose for $\psi$ a decomposition according to the support
\beq\label{E:psi-decomp}
\psi(k)=\eps^{1-d}\widehat{D}\left(\tfrac{k-k_0}{\eps}\right)+\eps^{1-d}\widehat{R}\left(\tfrac{k-k_0}{\eps}\right),
\eeq
where 
$$\text{supp}\widehat{D}\subset B_{\eps^{r-1}}:=\{k\in \R^d:|k|<\eps^{r-1}\}, \ \text{supp}\widehat{R}\subset B_{\eps^{r-1}}^c:=\R^d\setminus B_{\eps^{r-1}}$$
with $r\in (0,1)$. At the moment $r$ is a free parameter; it will be specified below. We also define 
$$\Bvechat_D(k) := \eps^{1-d}\widehat{D}\left(\tfrac{k-k_0}{\eps}\right) \vec{\eta}^{(j_0)}(k), \quad \Bvechat_R(k):= \eps^{1-d}\widehat{R}\left(\tfrac{k-k_0}{\eps}\right) \vec{\eta}^{(j_0)}(k)$$
such that
$$\Bvechat_P=\Bvechat_D+\Bvechat_R.$$
In the ansatz in \eqref{E:psi-decomp} we wish to find the component $\widehat{D}$ close to $\chi_{B_{\eps^{r-1}}} \widehat{C}$ and the component $\widehat{R}$ small. The component $\Bvechat_D$ is then approximated by $\Bvechat_\text{app}$ in \eqref{E:ans}. If also $\Bvechat_Q$ is small, then the whole constructed solution $\Bvechat$ is close to $\Bvechat_\text{app}$.

The Fourier transform of \eqref{E:CME} is
\beq\label{E:CME-FT}
\omega\Bvechat(k) - L(\ri k)\Bvechat(k) + \Nvechat(\Bvechat)(k)=0, \ k\in \R^d.
\eeq
For the selected ansatz equation \eqref{E:CME-FT} becomes
\begin{align}
\eps^{1-d}(\omega_0+\eps^2\omega_1-\lambda_{j_0}(k))\left(\widehat{D}\left(\tfrac{k-k_0}{\eps}\right) + \widehat{R}\left(\tfrac{k-k_0}{\eps}\right)\right)+\vec{\eta}^{(j_0)}(k)^*\Nvechat(\Bvechat)(k)&=0,\\
Q_k\left((\omega_0+\eps^2\omega_1)I-L(\ri k)\right)Q_k\Bvechat_Q(k) + Q_k\Nvechat(\Bvechat)(k)&=0.\label{E:Q-eq}
\end{align}
Because $\omega_0\in \sigma(L(\ri k_0))$, the inverse of the matrix $(\omega_0+\eps^2\omega_1)I-L(\ri k)$ is not bounded uniformly in $\eps$. In a neighbourhood of $k_0$ the norm of the inverse blows up as $\eps\to 0$. However, 
$$M_k:=Q_k\left((\omega_0+\eps^2\omega_1)I-L(\ri k)\right)Q_k$$ 
is invertible uniformly in $\eps$ due to assumption (A.2).

We separate the explicit part of \eqref{E:Q-eq} by writing
$$\Nvechat(\Bvechat)(k)=\Nvechat(\Bvechat_P)(k)+(\Nvechat(\Bvechat)(k)-\Nvechat(\Bvechat_P)(k))$$
and
$$\Bvechat_Q=\Bvechat_{Q,1}+\Bvechat_{Q,2},$$
where $\Bvechat_{Q,1}$ solves the explicit part, i.e.
\beq\label{E:BQ1-eq}
\Bvechat_{Q,1}(k) = -M_k^{-1} Q_k\Nvechat(\Bvechat_P)(k).
\eeq
The system to solve is thus
\begin{align}
\eps^{1-d}(\omega_0+\eps^2\omega_1-\lambda_{j_0}(k))\left(\widehat{D}\left(\tfrac{k-k_0}{\eps}\right) + \widehat{R}\left(\tfrac{k-k_0}{\eps}\right)\right)+\vec{\eta}^{(j_0)}(k)^*\Nvechat(\Bvechat)(k)&=0\label{E:P-eq}\\
M_k\Bvechat_{Q,2}(k) + Q_k(\Nvechat(\Bvechat)(k)-\Nvechat(\Bvechat_P)(k))&=0.\label{E:Qred-eq}
\end{align}
Our procedure for constructing a solution can be sketched as follows.
\ben
\item For any $\widehat{D},\widehat{R}\in L^1(\R^d)$ equation \eqref{E:BQ1-eq} produces a small $\Bvechat_{Q,1}$ (because $\|\Bvechat_P\|_{L^1}=O(\eps)$).
\item For any $\widehat{D},\widehat{R}\in L^1(\R^d)$ and $\Bvechat_{Q,1}$ from step 1 we solve \eqref{E:Qred-eq} by a fixed point argument for a small $\Bvechat_{Q,2}$.
\item For any $\widehat{D}\in L^1(\R^d)$ and for $\Bvechat_{Q}$ from steps 1 and 2 we solve \eqref{E:P-eq} with $k\in B_{\eps^{r}}(k_0)^c$ for a small $\widehat{R}$ by a fixed point argument.
\item With the components obtained in the above steps we find a solution $\widehat{D}\in L_2^1(\R^d)$ of \eqref{E:P-eq} with $k\in B_{\eps^{r}}(k_0)$ close to a $\widehat{C} \in L^1_2(\R^d)$ (with $C$ a solution of \eqref{E:NLS}) - provided such a $C$ exists. In addition $C$ needs to satisfy a certain symmetry, the $PT$-symmetry. Also here a fixed point argument is used  - roughly speaking for the difference $\widehat{D}-\widehat{C}$.
\item The error $\|\Bvechat - \Bvechat_\text{app}\|_{L^1(\R^d)}$ is $O(\eps^2)$ if $\widehat{C}$ decays fast enough, namely if $\widehat{C}\in L^1_4(\R^d)$.
\een

\blem\label{L:BP-est}
If $\widehat{D}\in L^1(\R^d)$, $\widehat{R}\in L^1_{s_R}(\R^d)$ for some $s_R\geq 0$, and $\text{supp}\widehat{D}\subset B_{\eps^{r-1}},$  $\text{supp}\widehat{R}\subset B_{\eps^{r-1}}^c$, then there are constant $c_1,c_2>0$ such that for all $\eps>0$ small enough
$$\|\Bvechat_P\|_{L^1(\R^d)}\leq c_1 \eps \left(\|\widehat{D}\|_{L^1(B_{\eps^{r-1}})}+\eps^{(1-r)s_R}\|\widehat{R}\|_{L_{s_R}^1(B_{\eps^{r-1}}^c)}\right)
$$
and
$$\|\Nvechat(\Bvechat_P)\|_{L^1(\R^d)}\leq c_2 \eps^3\left(\|\widehat{D}\|_{L^1(B_{\eps^{r-1}})}+\eps^{(1-r)s_R}\|\widehat{R}\|_{L_{s_R}^1(B_{\eps^{r-1}}^c)}\right)^3.$$
\elem
\bpf
Because $\|\widehat{f}(\eps^{-1}(\cdot-k_0))\|_{L^1(\R^d)}=\eps^d\|\widehat{f}\|_{L^1(\R^d)}$, we have
$$
\begin{aligned}
\|\Bvechat_P\|_{L^1(\R^d)} &\leq c\eps \left(\|\widehat{D}\|_{L^1(\R^d)} + \int_{\R^d}(1+|\kappa|)^{-s_R}(1+|\kappa|)^{s_R}|\widehat{R}(\kappa)|\dd\kappa\right)\\
&\leq c\eps \left(\|\widehat{D}\|_{L^1(\R^d)} + \sup_{\kappa\in B_{\eps^{r-1}}^c}(1+|\kappa|)^{-s_R}\|\widehat{R}\|_{L_{s_R}^1(B_{\eps^{r-1}}^c)}\right) \\
&\leq c_1 \eps \left(\|\widehat{D}\|_{L^1(B_{\eps^{r-1}})}+\eps^{(1-r)s_R}\|\widehat{R}\|_{L_{s_R}^1(B_{\eps^{r-1}}^c)}\right).
\end{aligned}
$$
The estimate for $\|\Nvechat(\Bvechat_P)\|_{L^1(\R^d)}$  follows by Young's inequality for convolutions.
\epf
\bi
\item[(1)] \underline{Component $\Bvechat_{Q,1}$} 
\ei
Because $\Bvechat_{Q,1}=-M_k^{-1}Q_k\Nvechat(\Bvechat_P)(k)$, we get from Lemma
\ref{L:BP-est} the estimate
\beq\label{E:BQ1-est}
\|\Bvechat_{Q,1}\|_{L^1(\R^d)}\leq c \eps^3 (\|\widehat{D}\|_{L^1(\R^d)}+\|\widehat{R}\|_{L^1(\R^d)})^3.
\eeq
\bi
\item[(2)] \underline{Component $\Bvechat_{Q,2}$} 
\ei
With $\Bvechat_{Q,1}$ from above (and $\widehat{D},\widehat{R}$ given) component $\Bvechat_{Q,2}$ satisfies
$$\Bvechat_{Q,2}(k) = M_k^{-1}Q_k\left( \Nvechat(\Bvechat_P)(k)-\Nvechat(\Bvechat_P+\Bvechat_{Q,1}+\Bvechat_{Q,2})(k) \right)=:\vec{G}(\Bvechat_{Q,2})(k).$$
Due to the cubic structure of $\vec{N}$ we have
$$\|\vec{G}(\Bvechat_{Q,2})\|_{L^1(\R^d)} \leq c \|\Bvechat_{Q,1}+\Bvechat_{Q,2}\|_{L^1(\R^d)}(\|\Bvechat_{P}\|_{L^1(\R^d)}+\|\Bvechat_{Q,1}\|_{L^1(\R^d)}+\|\Bvechat_{Q,2}\|_{L^1(\R^d)})^2.$$
For $\Bvechat_{Q,2} \in B_{\eps^\eta}, \eta >0$ we have
$$\|\vec{G}(\Bvechat_{Q,2})\|_{L^1(\R^d)} \leq c(\|\widehat{D}\|_{L^1},\|\widehat{R}\|_{L^1}) (\eps^5 + \eps^{2+\eta} + \eps^{1+2\eta}+ \eps^{3\eta}).$$
Hence, $\vec{G}:B^{(L^1)}_{c_0\eps^5}\to B^{(L^1)}_{c_0\eps^5}$ for some $c_0(\|\widehat{D}\|_{L^1},\|\widehat{R}\|_{L^1})>0$, where $B^{(L^1)}_\alpha := \{\vec{v}\in L^1(\R^d):\|\vec{v}\|_{L^1(\R^d)}\leq \alpha\}.$ The constants $c$ and $c_0$ depend polynomially on $\|\widehat{D}\|_{L^1}$ and $\|\widehat{R}\|_{L^1}$.

Similary, we obtain the contraction (for $\eps>0$ small enough)
$$\|\vec{G}(\Bvechat^{(1)}_{Q,R})-\vec{G}(\Bvechat^{(2)}_{Q,R})\|_{L^1(\R^d)} \leq c\eps^2\|\Bvechat^{(1)}_{Q,R}-\Bvechat^{(2)}_{Q,R}\|_{L^1(\R^d)},$$
if $\Bvechat^{(1)}_{Q,R},\Bvechat^{(2)}_{Q,R}\in B^{(L^1)}_{c_0\eps^5}$. Hence, for $\eps>0$ small enough we have a unique solution $\Bvechat_{Q,2}\in B^{(L^1)}_{c_0\eps^5}$ of $\Bvechat_{Q,2} =\vec{G}(\Bvechat_{Q,2})$, i.e.
\beq\label{E:BQR-est}
\|\Bvechat_{Q,2}\|_{L^1(\R^d)}\leq c_0(\|\widehat{D}\|_{L^1},\|\widehat{R}\|_{L^1})\eps^5.
\eeq
\bi
\item[(3)] \underline{Component $\Bvechat_R$} 
\ei
For any $\widehat{D}\in L^1(\R^d)$ and with the above estimate on $\Bvechat_{Q}$ we look for a small $\widehat{R}$. The support of $\widehat{R}$ is $B_{\eps^{r-1}}^c$, whence for $k\in \R^d\setminus B_{\eps^r}(k_0)$ we can divide in \eqref{E:P-eq} by $\omega_0+\eps^2\omega_1-\lambda_{j_0}(k)$ and obtain
$$\widehat{R}\left(\frac{k-k_0}{\eps}\right)= \eps^{d-1}\nu(k) \vec{\eta}^{(j_0)}(k)^*\Nvechat(\Bvechat)(k), \quad k \in \R^d\setminus B_{\eps^r}(k_0)$$
with
$$\nu(k):=(\lambda_{j_0}(k)-\omega_0-\eps^2\omega_1)^{-1}.$$ 
Since $\nabla \lambda_{j_0}(k_0)=0$, we have $|\lambda_{j_0}(k)-\omega_0|>c\eps^{2r}$ for all $k\in \R^d\setminus B_{\eps^r}(k_0)$ and hence $|\nu(k)|\leq c\eps^{-2r}$.
In order to exploit the localized nature of $\Bvechat_D$ and the smallness of $\Bvechat_Q$, we write for $\kappa:=\tfrac{k-k_0}{\eps}$
$$\widehat{R}(\kappa)=H(\widehat{R})(\kappa),$$
where 
\beq
\label{E:H}
\begin{aligned}
H(\widehat{R})(\kappa):=\eps^{d-1}\nu(k_0+\eps \kappa)&\left[h(k_0+\eps \kappa)^{-1}h(k_0+\eps \kappa)\vec{\eta}^{(j_0)}(k_0+\eps \kappa)^*\Nvechat(\Bvechat_D)(k_0+\eps \kappa)\right.\\
&\qquad \left. +\vec{\eta}^{(j_0)}(k_0+\eps \kappa)^*\left(\Nvechat(\Bvechat)(k_0+\eps \kappa)-\Nvechat(\Bvechat_D)(k_0+\eps \kappa)\right)\right]
\end{aligned}
\eeq
with $h(k):=\left(1+\tfrac{|k-k_0|}{\eps}\right)^{s_D}$. Using $\sup_{k\in \text{supp}\Bvechat_R}|h(k)^{-1}|\leq c\eps^{(1-r)s_D}$, the cubic form of $\Nvechat$ and the fact that $\left(\widehat{D}\left(\tfrac{\cdot -k_0}{\eps}\right)*\widehat{\overline{D}}\left(\tfrac{\cdot +k_0}{\eps}\right)*\widehat{D}\left(\tfrac{\cdot -k_0}{\eps}\right)\right)(k_0+\eps \kappa) = \eps^{2d}(\widehat{D}*\widehat{\overline{D}}*\widehat{D})(\kappa)$, we get
$$
\begin{aligned}
\|H(\widehat{R})\|_{L^1(\R^d)}\leq &c_1\eps^{2-2r+(1-r)s_D}\|(1+|\cdot|)^{s_D}(\widehat{D}*\widehat{\overline{D}}*\widehat{D})(\cdot)\|_{L^1(\R^d)}\\
 & + c_2\eps^{-2r-1}\left(\|\Bvechat_R\|_{L^1(\R^d)}+\|\Bvechat_Q\|_{L^1(\R^d)}\right)\left(\|\Bvechat_D\|_{L^1(\R^d)}+\|\Bvechat_R\|_{L^1(\R^d)}+\|\Bvechat_Q\|_{L^1(\R^d)}\right)^2.
\end{aligned}
$$
Since $\|\Bvechat_R\|_{L^1(\R^d)}\leq c \eps \|\widehat{R}\|_{L^1(\R^d)}$, $\|\Bvechat_D\|_{L^1(\R^d)}\leq c \eps \|\widehat{D}\|_{L^1(\R^d)}$, and $\|\Bvechat_Q\|_{L^1(\R^d)}\leq c(\|\widehat{D}\|_{L^1},\|\widehat{R}\|_{L^1})\eps^3$ (with a polynomial $c$), we have
$$
\begin{aligned}
\|H(\widehat{R})\|_{L^1(\R^d)}\leq &c\eps^{2-2r}\left(\eps^{(1-r)s_D}\|\widehat{D}\|^3_{L^1_{s_D}(\R^d)} + \|\widehat{R}\|_{L^1(\R^d)}^3+ \|\widehat{R}\|_{L^1(\R^d)}^2 + \|\widehat{R}\|_{L^1(\R^d)} + \eps^2\right)
\end{aligned}
$$
with $c$ depending polynomially on $\|\widehat{D}\|_{L^1}$ and $\|\widehat{R}\|_{L^1}$.

If $\widehat{D}\in L^1_{s_D}(\R^d)$, $s_D\geq 0$, then there is a constant $c$ depending polynomially on $\|\widehat{D}\|_{L^1_{s_D}(\R^d)}$ such that
\beq\label{E:H-and-alpha}
\|H(\widehat{R})\|_{L^1(\R^d)}\leq c(\|\widehat{D}\|_{L^1_{s_D}})\eps^\alpha, \quad \alpha:=\min\{4-2r,(2+s_D)(1-r)\}
\eeq
for all $R$ with $\|\widehat{R}\|_{L^1(\R^d)}\leq c\eps^\alpha$ and all $\eps>0$ small enough. Hence $H:B^{(L^1)}_{c\eps^\alpha}\to B^{(L^1)}_{c\eps^\alpha}$ for some $c>0$ if $\eps>0$ is small enough and if $\widehat{D}\in L^1_{s_D}(\R^d)$. 

Similarly, we get the contraction property of $H$ on $B^{(L^1)}_{c\eps^\alpha}$ for $\eps>0$ small enough. The constructed fixed point $\widehat{R}\in B^{(L^1)}_{c\eps^\alpha}$ yields for any $s_D\geq 0$
\beq\label{E:BR-est}
\|\Bvechat_R\|_{L^1(\R^d)} \leq c(\|\widehat{D}\|_{L^1_{s_D}})\eps^{\alpha+1}.
\eeq
\bi
\item[(4)] \underline{Component $\Bvechat_D$} 
\ei
Finally, we consider the component $\Bvechat_D$. For $k\in \text{supp}(\Bvechat_D)= B_{\eps^r}(k_0)$ we rewrite \eqref{E:P-eq} as follows. We add and subtract $\Nvechat(\Bvechat_D)$ like in \eqref{E:H}, we Taylor expand $\lambda_{j_0}(k)$ at $k=k_0$, and we use the variable $\kappa=\eps^{-1}(k-k_0)$. This leads to
\beq\label{E:D-eq}
(\omega_1-\kappa^T G_0\kappa)\widehat{D}(\kappa)+\chi_{B_{\eps^{r-1}}}(\kappa)\Gamma (\widehat{D}*\widehat{\overline{D}}*\widehat{D})(\kappa)=\rho(\widehat{D})(\kappa),
\eeq
where
$$
\begin{aligned}
&\rho(\widehat{D})(\kappa) = \chi_{B_{\eps^{r-1}}}(\kappa)\left[\Gamma (\widehat{D}*\widehat{\overline{D}}*\widehat{D})(\kappa) -\eps^{d-3} \vec{\eta}^{(j_0)}(k_0+\eps \kappa)^*\Nvechat(\Bvechat_D)(k_0+\eps \kappa)\right]\\
& +\eps^{-2}\left[\delta(k_0+\eps \kappa)\widehat{D}(\kappa) + \eps^{d-1}\chi_{B_{\eps^{r-1}}}(\kappa) \vec{\eta}^{(j_0)}(k_0+\eps \kappa)^*\left(\Nvechat(\Bvechat)(k_0+\eps \kappa)-\Nvechat(\Bvechat_D)(k_0+\eps \kappa)\right)\right]
\end{aligned}
$$
and $\delta(k):=\lambda_{j_0}(k)-\omega_0-\frac{1}{2}(k-k_0)^TD^2\lambda_{j_0}(k_0)(k-k_0)$.

Next, we write $\rho=\rho_1+\rho_2+\rho_3$, where 
$$
\begin{aligned}
\rho_1(\widehat{D})(\kappa)&:=\chi_{B_{\eps^{r-1}}}(\kappa)\left[\Gamma (\widehat{D}*\widehat{\overline{D}}*\widehat{D})(\kappa) -\eps^{d-3} \vec{\eta}^{(j_0)}(k_0+\eps \kappa)^*\Nvechat(\Bvechat_D)(k_0+\eps \kappa)\right],\\
\rho_2(\widehat{D})(\kappa)&:= \eps^{-2} \delta(k_0+\eps \kappa)\widehat{D}(\kappa),\\
\rho_3(\widehat{D})(\kappa)&:= \eps^{d-3}\chi_{B_{\eps^{r-1}}}(\kappa) \vec{\eta}^{(j_0)}(k_0+\eps \kappa)^*\left(\Nvechat(\Bvechat)(k_0+\eps \kappa)-\Nvechat(\Bvechat_D)(k_0+\eps \kappa)\right).
\end{aligned}
$$
For $\rho_1$ we get
$$
\begin{aligned}
\rho_1(\widehat{D})(\kappa)=&\eps^{-2d}\sum_{\stackrel{m,n,o\in\{1,\dots,N\}}{j\in\{1,\dots,N\}}}\gamma_j^{(m,n,o)}\int_{B_{2\eps^r}}\int_{B_{\eps^r}(k_0)}\left(\overline{\eta}_j^{(j_0)}(k_0+\eps\kappa)\eta_m^{(j_0)}(k_0+\eps\kappa-s)\overline{\eta}_n^{(j_0)}(s-t)\eta_o^{(j_0)}(t)\right.\\
&\left.-\Gamma\right)\widehat{D}\left(\kappa-\frac{s}{\eps}\right)\widehat{\overline{D}}\left(\frac{s-t+k_0}{\eps}\right)\widehat{D}\left(\frac{t-k_0}{\eps}\right)\dd t  \dd s\\
=&\sum_{\stackrel{m,n,o\in\{1,\dots,N\}}{j\in\{1,\dots,N\}}}\gamma_j^{(m,n,o)}\int_{B_{2\eps^{r-1}}}\int_{B_{\eps^{r-1}}}\left(\overline{\eta}_j^{(j_0)}(k_0+\eps\kappa)\eta_m^{(j_0)}(k_0+\eps(\kappa-\tilde{s}))\overline{\eta}_n^{(j_0)}(k_0+\eps(\tilde{s}-\tilde{t}))\times\right.\\
&\left.\times \eta_o^{(j_0)}(k_0+\eps\tilde{t}) -\Gamma\right)\widehat{D}\left(\kappa-\tilde{s}\right)\widehat{\overline{D}}\left(\tilde{s}-\tilde{t}\right)\widehat{D}\left(\tilde{t}\right)\dd \tilde{t}  \dd \tilde{s}.
\end{aligned}
$$
Because $\Gamma=\sum_{m,n,o,j\in\{1,\dots,N\}}\gamma_j^{(m,n,o)}\overline{\eta}_j^{(j_0)}(k_0)\eta_m^{(j_0)}(k_0)\overline{\eta}_n^{(j_0)}(k_0)\eta_o^{(j_0)}(k_0)$, we get 
$$
\begin{aligned}
&\left|\sum_{m,n,o,j\in\{1,\dots,N\}}\gamma_j^{(m,n,o)}\overline{\eta}_j^{(j_0)}(k_0+\eps\kappa)\eta_m^{(j_0)}(k_0+\eps(\kappa-\tilde{s}))\overline{\eta}_n^{(j_0)}(k_0+\eps(\tilde{s}-\tilde{t}))\eta_o^{(j_0)}(k_0+\eps\tilde{t}) -\Gamma\right|\\
&\leq c\eps(|\kappa|+|\kappa-\tilde{s}|+|\tilde{s}-\tilde{t}|+|\tilde{t}|)
\end{aligned}
$$
due to the Lipschitz continuity of $k\mapsto\eta^{(j_0)}(k)$. Hence, by Young's inequality for convolutions,
\beq\label{E:rho1-est}
\|\rho_1(\widehat{D})\|_{L^1(\R^d)} \leq c\eps\|\widehat{D}\|_{L^1_1(\R^d)}^3.
\eeq

For $\rho_2$ we note that $|\delta(k_0+\eps \kappa)|\leq c\eps^3|\kappa|^3$ for $\kappa\in B_{\eps^{r-1}}$ and $\eps>0$ small enough. We estimate
\begin{align}
\|\rho_2(\widehat{D})\|_{L^1(\R^d)} &\leq c\eps \int_{B_{\eps^{r-1}}}|\kappa|^3|\widehat{D}(\kappa)|\dd \kappa\leq c\eps\sup_{\kappa\in B_{\eps^{r-1}}}|\kappa|^\beta \int_{\R^d}|\kappa|^{3-\beta}|\widehat{D}(\kappa)|\dd \kappa \notag\\
&\leq c\eps^{1-\beta(1-r)}\|\widehat{D}\|_{L^1_{3-\beta}(\R^d)}\label{E:rho2-est}
\end{align}
for any $\beta\in[0,3)$.

Finally, we estimate $\rho_3$. Note that $\Nvechat(\Bvechat)(k_0+\eps \cdot)-\Nvechat(\Bvechat_D)(k_0+\eps \cdot)$ appears also in \eqref{E:H}. We have
$$
\begin{aligned}
&\|\Nvechat(\Bvechat)(k_0+\eps \cdot)-\Nvechat(\Bvechat_D)(k_0+\eps \cdot)\|_{L^1(\R^d)}\\
&\leq \eps^{-d}\left(\|\Bvechat_R\|_{L^1(\R^d)}+\|\Bvechat_Q\|_{L^1(\R^d)}\right)\left(\|\Bvechat_D\|_{L^1(\R^d)}+\|\Bvechat_R\|_{L^1(\R^d)}+\|\Bvechat_Q\|_{L^1(\R^d)}\right)^2\\
&\leq c_1(\|\widehat{D}\|_{L^1_{s_D}})\eps^{-d}\left((\eps^{\alpha+3}+\eps^5)\|\widehat{D}\|^2_{L^1(\R^d)}+(\eps^{2\alpha+3}+\eps^7)\|\widehat{D}\|_{L^1(\R^d)}+\eps^{3\alpha+3}+\eps^9\right) \\
&\leq c_2(\|\widehat{D}\|_{L^1_{s_D}})(\eps^{\alpha+3-d}+\eps^{5-d})
\end{aligned}
$$
for $\eps>0$ small enough, where we have made use of \eqref{E:BQ1-est}, \eqref{E:BQR-est}, \eqref{E:BR-est}, the fact that $\|\Bvechat_D\|_{L^1}\leq \eps \|\widehat{D}\|_{L^1}$ and the estimate $\|\widehat{R}\|_{L^1}\leq c(\|\widehat{D}\|_{L^1_{s_D}})\eps^\alpha$. The dependence of $c_1$ and $c_2$ on $\|\widehat{D}\|_{L^1_{s_D}}$ is polynomial. As a result
\beq\label{E:rho3-est}
\|\rho_3(\widehat{D})\|_{L^1(\R^d)} \leq c(\|\widehat{D}\|_{L^1_{s_D}})(\eps^{\alpha}+\eps^{2}).
\eeq

The whole right hand side of \eqref{E:D-eq} is thus estimated as
$$\|\rho(\widehat{D})\|_{L^1(\R^d)} \leq c(\|\widehat{D}\|_{L^1_{s_D}},\|\widehat{D}\|_{L^1_{3-\beta}}) \left(\eps^{1-\beta(1-r)}+\eps^{\alpha}\right)$$ 
for any $\beta \in [0,3)$. Once again, the constant $c$ depends polynomially on its arguments.

In order to solve \eqref{E:D-eq} for $\widehat{D}$ below (using a fixed point argument), we need to consider $\rho$ as an inhomogeneity. The linearized operator to be inverted in the iteration is of second order such that in Fourier space it acts from $L^1_2(\R^d)$ to $L^1(\R^d)$. For that reason we need to choose $s_D=2$ and $\beta \geq 1$ above. The choice $s_D=2,  \beta=1, r=1/2$ leads to $\alpha=2$ and $1-\beta(1-r)=1/2$, i.e.
$\|\rho(\widehat{D})\|_{L^1(\R^d)} \leq c \eps^{1/2}$ provided $\widehat{D}\in L^1_{2}(\R^d)$. For $s_D=2,\beta =1$ the largest value of $ \min\{1-\beta(1-r),\alpha\}$ is $4/5$ attained at $r=4/5$. Hence,  with $r=4/5$ we get  the best possible estimate
$$
\|\rho(\widehat{D})\|_{L^1(\R^d)} \leq c(\|\widehat{D}\|_{L^1_{2}}) \eps^{4/5}.
$$
This order determines the accuracy of the approximation and turns out to be insufficient. It leads to $
\|\widehat{\vec{B}}-\Bvechat_{\text{app}}\|_{L^1(\R^d)}\leq c\eps^{9/5}$ instead of $c\eps^2$.

Clearly, the leading order term in the residual $\rho$ is caused by the error from the Taylor expansion of $\lambda_{j_0}$. We introduce a refined ansatz for $\widehat{D}$ in order to make this error of higher order. Note that equation \eqref{E:D-eq} is a perturbation of the NLS \eqref{E:NLS}. Writing 
$$f_\text{NLS}(C):=\omega_1 C + \nabla^T(G_0\nabla C) + \Gamma |C|^2C,$$
equation \eqref{E:D-eq} is
\beq\label{E:Dhat-eq}
\chi_{B_{\eps^{r-1}}}\widehat{f_\text{NLS}(D)} = \rho(\widehat{D}).
\eeq
We search for $D$ close to a solution $C$ of the NLS, i.e. of $f_\text{NLS}(C)=0$. For that we need to look for $\widehat{D}$ in a vicinity of $$\widehat{C}^{(\eps)}:=\chi_{B_{\eps^{r-1}}}\widehat{C}.$$ 
We choose the following ansatz
\beq\label{Dhat-ans-improve}
\widehat{D}(\kappa)=\left(1+\eps\frac{D^3\lambda_{j_0}(k_0)(\kappa,\kappa,\kappa)}{6(\omega_1-\kappa^TG_0\kappa)}\right)\widehat{C}^{(\eps)}(\kappa)+\widehat{d}(\kappa),
\eeq
where $D^3\lambda_{j_0}(k_0)(\kappa,\kappa,\kappa)$ is the third order term in the Taylor expansion of $\eps^{-3}\lambda_{j_0}(k_0+\eps\kappa)$ and where $\text{supp}(\widehat{d}) \subset B_{\eps^{r-1}}$. We look for a solution with a small $\widehat{d}$.

We also define 
$$\widehat{D}_0:=\widehat{C}^{(\eps)}+\widehat{d}, \quad \nu(\kappa):=\frac{D^3\lambda_{j_0}(k_0)(\kappa,\kappa,\kappa)}{6(\omega_1-\kappa^TG_0\kappa)},$$
such that $\widehat{D}=\widehat{D}_0+\eps \nu \widehat{C}^{(\eps)}$. 
With this notation equation \eqref{E:Dhat-eq} reads
\beq\label{E:D0hat-eq}
\chi_{B_{\eps^{r-1}}}\widehat{f_\text{NLS}(D_0)} = \rho(\widehat{D}) + \chi_{B_{\eps^{r-1}}}\left(\widehat{f_\text{NLS}(D_0)}-\widehat{f_\text{NLS}(D)}\right)=:\tilde{\rho}(\widehat{D}).
\eeq
Compared to $\rho$ the right hand side $\tilde{\rho}$ is smaller as we show next. It is
$$
\tilde{\rho}(\widehat{D})=\rho_1(\widehat{D})+\rho_3(\widehat{D})+\tilde{\rho}_2(\widehat{D}),
$$
where 
$$
\begin{aligned}
\tilde{\rho}_2(\widehat{D})(\kappa)=& \eps^{-2} \left(\lambda_{j_0}(k_0+\eps \kappa)-\omega_0-\eps^2\kappa^TG_0\kappa\right)\widehat{D}(\kappa)-\eps(\omega_1-\kappa^TG_0\kappa)\nu(\kappa)\widehat{C}^{(\eps)}(\kappa)\\
&+ \chi_{B_{\eps^{r-1}}}(\kappa)\Gamma \left[\widehat{D}_0*\widehat{\overline{D}_0}*\widehat{D}_0-\widehat{D}*\widehat{\overline{D}}*\widehat{D}\right](\kappa)\\
=&  \eps^{-2} \left(\lambda_{j_0}(k_0+\eps \kappa)-\omega_0-\eps^2\kappa^TG_0\kappa -\frac{\eps^3}{6} D^3\lambda_{j_0}(k_0)(\kappa,\kappa,\kappa)\right)\widehat{C}^{(\eps)}(\kappa)\\
&+\eps^{-2} \left(\lambda_{j_0}(k_0+\eps \kappa)-\omega_0-\eps^2\kappa^TG_0\kappa\right)\left(\eps\nu(\kappa)\widehat{C}^{(\eps)}(\kappa)+ \widehat{d}(\kappa)\right)\\
&+ \chi_{B_{\eps^{r-1}}}(\kappa)\Gamma \left[\widehat{D}_0*\widehat{\overline{D}_0}*\widehat{D}_0-\widehat{D}*\widehat{\overline{D}}*\widehat{D}\right](\kappa).
\end{aligned}
$$
To estimate $\tilde{\rho}_2$ note that the first Taylor expansion error term (i.e. the first line) can be estimated in $L^1(\R^d)$ by $c\eps^2\|\widehat{C}^{(\eps)}\|_{L^1_4(\R^d)}\leq c\eps^2\|\widehat{C}\|_{L^1_4(\R^d)}$. For the second Taylor expansion error note first that $|\nu(\kappa)|\leq c(1+|\kappa|)$ for all $\kappa\in \R^d$ such that
$$\|\eps^{-2} \left(\lambda_{j_0}(k_0+\eps \cdot)-\omega_0-\eps^2(\cdot)^TG_0(\cdot)\right)(\eps\nu\widehat{C}^{(\eps)}+ \widehat{d})\|_{L^1(\R^d)}\leq c\left(\eps^2\|\widehat{C}\|_{L^1_4(\R^d)} + \eps \int_{\R^d}|\kappa|^3|\widehat{d}(\kappa)|\dd\kappa\right).$$
Just like in \eqref{E:rho2-est} we get 
$$\eps \int_{\R^d}|\kappa|^3|\widehat{d}(\kappa)|\dd\kappa \leq  c\eps^{1-\beta(1-r)}\|\widehat{d}\|_{L^1_{3-\beta}(\R^d)}$$
for any $\beta \in [0,3)$. Once again, we need to choose $\beta \geq 1$ such that $\tilde{\rho}_2:\widehat{d}\mapsto \tilde{\rho}_2(\widehat{d})$ is a mapping from $L^1_2(\R^d)$ to $L^1(\R^d)$. With $\beta =1$ and $r=1/2$ we get
$1-\beta(1-r)=1/2$.

The last term in $\tilde{\rho}_2$ is 
$$
\begin{aligned}
&\chi_{B_{\eps^{r-1}}}\Gamma\left[\widehat{D}_0*\widehat{\overline{D}_0}*\widehat{D}_0-\widehat{D}*\widehat{\overline{D}}*\widehat{D}\right] = \\
&\chi_{B_{\eps^{r-1}}}\Gamma\left[(\widehat{D}-\eps \nu \widehat{C}^{(\eps)})*(\widehat{\overline{D}}-\eps \nu \widehat{\overline{C}}^{(\eps)})*(\widehat{D}-\eps \nu \widehat{C}^{(\eps)})-\widehat{D}*\widehat{\overline{D}}*\widehat{D}\right]. 
\end{aligned}
$$
In the $L^1$-norm this can be estimated using Young's inequality by
$$c\left(\eps\|\widehat{C}\|_{L^1_1(\R^d)}\|\widehat{D}\|^2_{L^1(\R^d)}+\eps^2\|\widehat{C}\|^2_{L^1_1(\R^d)}\|\widehat{D}\|_{L^1(\R^d)}+\eps^3\|\widehat{C}\|^3_{L^1_1(\R^d)}\right).$$
In summary, with the above choice of $\beta$ and $r$ we get
\beq\label{E:rho2til-est}
\|\tilde{\rho}_2\|_{L^1(\R^d)}\leq c\left(\eps^{1/2}\|\widehat{d}\|_{L^1_2}+\eps^2 \|\widehat{C}\|_{L^1_4(\R^d)} +\eps\|\widehat{C}\|_{L^1_1(\R^d)}\|\widehat{D}\|^2_{L^1(\R^d)}+\eps^2\|\widehat{C}\|^2_{L^1_1(\R^d)}\|\widehat{D}\|_{L^1(\R^d)}+\eps^3\|\widehat{C}\|^3_{L^1_1(\R^d)}\right).\eeq
The terms $\rho_1$ and $\rho_3$ are estimated in \eqref{E:rho1-est} and \eqref{E:rho3-est}. As discussed above, we seek $\widehat{D}$ in $L^1_2(\R^d)$ and hence we set $s_D=2$ in \eqref{E:H-and-alpha}. This yields $\alpha=2$ and 
$$\|\rho_3(\widehat{D})\|_{L^1(\R^d)} \leq c(\|\widehat{D}\|_{L^1_2})\eps^2.$$
In summary, for any $\widehat{C}\in L^1_4(\R^d)$ fixed (with $C$ being a solution of the NLS), the right hand side of \eqref{E:D0hat-eq} is estimated as
\beq\label{E:est-rho-optim}
\|\tilde{\rho}(\widehat{D})\|_{L^1(\R^d)}\leq c \left(\eps + \|\widehat{d}\|_{L^1_2}\eps^{1/2} + (\|\widehat{d}\|_{L^1_2}^2+\|\widehat{d}\|_{L^1_2}^3)\eps\right).
\eeq

We proceed with a fixed point argument for the correction $d$. In order to obtain a differentiable function (to use the Jacobian of the NLS), we write $f_\text{NLS}$ in real variables. Writing $D=D_R+\ri D_I, C=C_R+\ri C_I$, $d=d_R+\ri d_I$, and $\tilde{\rho} = \tilde{\rho}_R+\ri \tilde{\rho}_I$, we define 
$$F_\text{NLS}(D_R,D_I):=\begin{pmatrix}\Real(f_\text{NLS}(D_R+\ri D_I))\\\Imag(f_\text{NLS}(D_R+\ri D_I))\end{pmatrix},$$
the Jacobian
\beq\label{E:J}
J:=DF_\text{NLS}(C_R,C_I)
\eeq
as well as the Fourier-truncation of the Jacobian
$$\widehat{J}_\eps:=\chi_{B_{\eps^{r-1}}} \cF\left(DF_\text{NLS}(C_R^{(\eps)},C_I^{(\eps)})\right),$$
where $\widehat{C_R}^{(\eps)}:=\chi_{B_{\eps^{r-1}}} \widehat{C_R}, C_R^{(\eps)}:=\cF^{-1}(\widehat{C_R}^{(\eps)})$ and $\widehat{C_I}^{(\eps)}:=\chi_{B_{\eps^{r-1}}} \widehat{C_I}$, $C_I^{(\eps)}:=\cF^{-1}(\widehat{C_I}^{(\eps)})$.
$\widehat{J}_\eps$ has the form
$$
\begin{aligned}
\widehat{J}_\eps = & \chi_{B_{\eps^{r-1}}} (\omega_1 +\kappa^T G_0 \kappa) \text{Id}_{2\times 2} \\
& + \chi_{B_{\eps^{r-1}}} \Gamma\begin{pmatrix} 3\widehat{C_R}^{(\eps)}*\widehat{C_R}^{(\eps)}* + \widehat{C_I}^{(\eps)}*\widehat{C_I}^{(\eps)}*    & 2\widehat{C_R}^{(\eps)}*\widehat{C_I}^{(\eps)}* \\
2\widehat{C_R}^{(\eps)}*\widehat{C_I}^{(\eps)}* & 3\widehat{C_I}^{(\eps)}*\widehat{C_I}^{(\eps)}* + \widehat{C_R}^{(\eps)}*\widehat{C_R}^{(\eps)}*     \end{pmatrix}.
\end{aligned}
$$

With this notation \eqref{E:D0hat-eq} reads
$$\widehat{J}_\eps \dvechat =W(\dvechat),$$
where
$$\dvechat:=(\widehat{d_R}, \widehat{d_I\hspace{0.05cm}})^T \text{ and } W(\dvechat):=\chi_{B_{\eps^{r-1}}}\vec{\tilde{\rho}}((1+\eps\nu)\widehat{\vec{C}}^{(\eps)}+\dvechat)- \chi_{B_{\eps^{r-1}}}\left(F_\text{NLS}(\widehat{\vec{C}}^{(\eps)}+\dvechat)-\widehat{J}_\eps \dvechat\right)$$
with $\widehat{\vec{C}}^{(\eps)}:=(\widehat{C_R}^{(\eps)},\widehat{C_I}^{(\eps)})^T$ and $\vec{\tilde{\rho}}:=(\tilde{\rho}_R,\tilde{\rho}_I)^T$.

The aim is to construct a small fixed point $\dvechat\in L^1_{s_D}(\R^d)$ of $\widehat{J}_\eps^{-1}W$. The difficulty is that the inverse of $\widehat{J}_\eps$ is not bounded uniformly in $\eps$. This is due to the presence of the $d+1$ zero eigenvalues of $J$ caused
by the $d$ spatial shift invariances and the phase invariance of the NLS, see assumption (A.3). The distance of the essential spectrum of $\widehat{J}_\eps$ from zero is $|\omega_1|$ due to the choice of $\text{sign}(\omega_1)$. To eliminate the zero eigenvalues, we work in a symmetric subspace of $L^1_{s_D}(\R^d)$ in which the invariances do not hold. A natural symmetry is the $PT$-symmetry. Hence, we consider the fixed point problem
$$ \dvechat =\widehat{J}_\eps^{-1} W(\dvechat)$$
in the space
$$
\begin{aligned}
X_{s_D}^\text{sym}:=&\{\dvechat\in L^1_{s_D}(\R^d): \text{supp}(\dvechat)\subset B_{\eps^{r-1}}, d(-x)=\overline{d(x)} \ \forall x \in \R^d\}\\
=& \{\dvechat\in L^1_{s_D}(\R^d): \text{supp}(\dvechat)\subset B_{\eps^{r-1}}, \Imag(\widehat{d_R})=-\Real(\widehat{d_I\hspace{0.05cm}})\}.
\end{aligned}
$$
Note that $\widehat{J}_0,\widehat{J}_\eps:L^1_{q}(\R^d)\to L^1_{q-2}(\R^d)$ for any $q\geq 2$. Because of assumption (A.3) $\widehat{J}_0^{-1}$ is bounded in $X_{s_D}^\text{sym}$ for any $s_D\geq 2$. Since $\widehat{J}_\eps$ is a perturbation of $\widehat{J}_0$, we still need to ensure that $0$ is not an eigenvalue of $\widehat{J}_\eps$. For that we use a spectral stability result of Kato, see \cite[Theorem IV.3.17]{Kato_1995}. Applied to our problem in the Banach space $X_{s_D-2}^\text{sym}, s_D\geq 2$ with the domain of  $\widehat{J}_0$  being $D(\widehat{J}_0)=X_{s_D}^\text{sym}$, it reads:

\textsl{Assume $\widehat{J}_\eps-\widehat{J}_0$ is $\widehat{J}_0$-bounded, i.e. $\mbox{domain}(\widehat{J}_0)\subset \mbox{domain}(\widehat{J}_\eps-\widehat{J}_0)$ and for some $a,b\geq 0$ is
\beq\label{E:Kato-cond1}
\|(\widehat{J}_\eps-\widehat{J}_0)\dvechat\|_{L^1_{s_D-2}} \leq a\|\dvechat\|_{L^1_{s_D-2}}+b\|\widehat{J}_0\dvechat\|_{L^1_{s_D-2}} \text{ for all } \dvechat\in X_{s_D}^\text{sym}.
\eeq
If for some $\zeta\in \rho (\widehat{J}_0)$
\beq\label{E:Kato-cond2}
a\|(\widehat{J}_0-\zeta)^{-1}\|_{L^1_{s_D-2}\to L^1_{s_D}}+b\|\widehat{J}_0(\widehat{J}_0-\zeta)^{-1}\|_{L^1_{s_D-2}\to L^1_{s_D-2}}<1,
\eeq
then $\zeta\in \rho (\widehat{J}_\eps)$ and
$$
\begin{aligned}
\|(\widehat{J}_\eps-\zeta)^{-1}\|_{L^1_{s_D-2}\to L^1_{s_D}}\leq &\|(\widehat{J}_0-\zeta)^{-1}\|_{L^1_{s_D-2}\to L^1_{s_D}}\left(1-a\|(\widehat{J}_0-\zeta)^{-1}\|_{L^1_{s_D-2}\to L^1_{s_D}}\right.\\
&\left.-b\|\widehat{J}_0(\widehat{J}_0-\zeta)^{-1}\|_{L^1_{s_D-2}\to L^1_{s_D-2}}\right)^{-1}.
\end{aligned}
$$}

We check now \eqref{E:Kato-cond1} and \eqref{E:Kato-cond2} for $\zeta =0$. For $\dvechat\in X_{s_D}^\text{sym}$ one has
$$
\begin{aligned}
\|(\widehat{J}_\eps-\widehat{J}_0)\dvechat\|_{L^1_{s_D-2}} = |\Gamma|&\left\|\left(\widehat{C_R}^{(\eps)}*\widehat{C_R}^{(\eps)}-\widehat{C_R}*\widehat{C_R}\right)*(3\widehat{d_R}+\widehat{d_I\hspace{0.05cm}})\right.\\
&+\left(\widehat{C_I}^{(\eps)}*\widehat{C_I}^{(\eps)}-\widehat{C_I}*\widehat{C_I}\right)*(3\widehat{d_I\hspace{0.05cm}}+\widehat{d_R})\\
&\left. +2\left(\widehat{C_R}^{(\eps)}*\widehat{C_I}^{(\eps)}-\widehat{C_R}*\widehat{C_I}\right)*(\widehat{d_R}+\widehat{d_I\hspace{0.05cm}})\right\|_{L^1_{s_D-2}}. 
\end{aligned}
$$
Writing $\widehat{C_{R,I}}^{(\eps)}=:\widehat{C_{R,I}}+\widehat{\gamma_{R,I}}^{(\eps)}$, it is $\text{supp}(\widehat{\gamma_{R,I}}^{(\eps)})\subset B_{\eps^{r-1}}^c$. The difference $(\widehat{J}_\eps-\widehat{J}_0)\dvechat$ consists of terms that are 
linear or quadratic in $\widehat{\gamma_{R,I}}^{(\eps)}$; for instance terms like $\widehat{\gamma_{R}}^{(\eps)}*\widehat{C_{R}}^{(\eps)}*\widehat{d_{R}}$ or $\widehat{\gamma_{I}}^{(\eps)}*\widehat{C_{R}}^{(\eps)}*\widehat{d_I\hspace{0.05cm}}$. Because
$$
\begin{aligned}
\|\widehat{\gamma_{R,I}}^{(\eps)}\|_{L^1_{s_D-2}}&=\int_{B_{\eps^{r-1}}^c} (1+|\kappa|)^{s_D-2}|\widehat{\gamma_{R,I}}(\kappa)|\dd\kappa\\
&\leq \sup_{|\kappa|>\eps^{r-1}}\frac{1}{(1+|\kappa|)^2}\|\widehat{\gamma_{R,I}}^{(\eps)}\|_{L^1_{s_D}}\leq c\eps^{2(1-r)}\|\widehat{\gamma_{R,I}}^{(\eps)}\|_{L^1_{s_D}},
\end{aligned}
$$
Young's inequality for convolutions yields
$$\|(\widehat{J}_\eps-\widehat{J}_0)\dvechat\|_{L^1_{s_D-2}} \leq c(\|\widehat{C}\|_{L^1_{s_D}})\eps^{2(1-r)}\|\dvechat\|_{L^1_{s_D-2}}$$
with $c$ depending polynomially on $\|\widehat{C}\|_{L^1_{s_D}}$. Conditions \eqref{E:Kato-cond1} and \eqref{E:Kato-cond2} for $\zeta =0$ are thus satisfied with $a=c\eps^{2(1-r)}$ and $b=0$ if $\eps>0$ is small enough and if $\widehat{C}\in L^1_{s_D}(\R^d)$.

For the fixed point problem we use $s_D=2$ and $r=1/2$ and show firstly that if $\widehat{C} \in X_{4}^\text{sym}$, then there is some $c>0$ such that for all $\eps>0$ small enough
\beq\label{E:selfmap}
\widehat{J}_\eps^{-1}W:B_{c\eps}^{2,\text{ sym}}\to B_{c\eps}^{2,\text{ sym}},
\eeq
where $B_{c\eps}^{2,\text{ sym}}$ is the $c\eps$-ball in $X_{2}^\text{sym}$, i.e.
$$B_{c\eps}^{2,\text{ sym}}:=\{\phi \in X_{2}^\text{sym}: \|\phi\|_{L^1_{2}}\leq c\eps\}.$$
Note that the requirement $\widehat{C}\in L^1_4$ is dictated by \eqref{E:rho2til-est}.

Secondly, we prove that $\widehat{J}_\eps^{-1}W$ is contractive provided $\widehat{C}\in \{\widehat{f}\in L^1_{4}(\R^d): \Imag(\widehat{f_R})=-\Real(\widehat{f_I}) \}$. We start by showing $\widehat{J}_\eps:X_{2}^\text{sym}\to X_{0}^\text{sym}$. The loss of $2$ in the weight is due to the second order nature of the operator $J$, i.e. due to the factor $\kappa^TG_0\kappa$ in Fourier variables. The entries $\omega_1+\kappa^TG_0\kappa$ clearly preserve the $PT$-symmetry. For the convolution terms we have, for instance
$$\Imag(\widehat{C_R})=-\Real(\widehat{C_I}) \ \Rightarrow \ \Imag(\widehat{C_R}^{(\eps)})=-\Real(\widehat{C_I}^{(\eps)}) \ \Rightarrow \ C^{(\eps)}(-x)=\overline{C^{(\eps)}(x)} \ \forall x.$$
Hence $C_R^{(\eps)^2}d_R, C_I^{(\eps)^2}d_R$, and $C_R^{(\eps)}C_I^{(\eps)}d_I$ are even and $C_I^{(\eps)^2}d_I, C_R^{(\eps)^2}d_I$, and $C_R^{(\eps)}C_I^{(\eps)}d_R$  are odd such that the $PT-$symmetry is preserved also by the convolution terms. In end effect, $\Imag ((\widehat{J}_\eps \dvechat)_1)=-\Real((\widehat{J}_\eps \dvechat)_2)$.
Hence, $\widehat{J}_\eps:X_{2}^\text{sym}\to X_{0}^\text{sym}$ and for $\widehat{J}_\eps^{-1}$ we get $\widehat{J}_\eps^{-1}:X_{0}^\text{sym}\to X_{2}^\text{sym}$.

Next, we show that $W:B_{c\eps}^{2,\text{ sym}}\to B_{c\eps}^{0,\text{sym}}$ if $\widehat{C} \in X_{4}^\text{sym}$. The term $\rho$ is estimated in \eqref{E:est-rho-optim} and dictates the order $\eps^1$. 

The difference $F_\text{NLS}(\widehat{\vec{C}}^{(\eps)}+\dvechat)-\widehat{J}_\eps \dvechat$ consists of terms quadratic in $\dvechat$  and hence is bounded in $L^1(\R^d)$ by $c_1(\|\dvechat\|_{L^1(\R^d)}^2+\|\dvechat\|_{L^1(\R^d)}^3)$. In summary,
$$\|W(\dvechat)\|_{L^1(\R^d)} \leq c_2 (\eps + \eps^{1/2}\|\dvechat\|_{L^1_2(\R^d)}
 + \|\dvechat\|_{L^1_2(\R^d)}^2+\|\dvechat\|_{L^1_2(\R^d)}^3) \leq  c \eps $$
if $\|\dvechat\|_{L^1_2(\R^d)} \leq c\eps$ and $\eps>0$ is small enough. Due to the boundedness of $\widehat{J}_\eps^{-1}: X_{0}^\text{sym}\to X_{2}^\text{sym}$ we thus have \eqref{E:selfmap}.

The contractive property of $\widehat{J}_\eps^{-1}W$ in $B_{c\eps}^{2,\text{ sym}}$ is now clear due to the quadratic nature of $F_\text{NLS}(\widehat{\vec{C}}^{(\eps)}+\dvechat)-\widehat{J}_\eps \dvechat$.

We conclude that if the solution $C$ of the NLS \eqref{E:NLS} satisfies $\widehat{C}\in \{\widehat{f}\in L^1_{4}(\R^d): \Imag(\widehat{C_R})=-\Real(\widehat{C_I}) \}$, then there is $c>0$ such that for all $\eps>0$ small enough the constructed solution $D$ of \eqref{E:D-eq} satisfies 
$$\widehat{D} \in \{\widehat{f}\in L^1_{2}(\R^d): \Imag(\widehat{f_R})=-\Real(\widehat{f_I\hspace{0.05cm}}) \}$$
and due to \eqref{Dhat-ans-improve}
\beq\label{E:D-C-est}
\|\widehat{D}-\chi_{B_{\eps^{-1/2}}}\widehat{C}\|_{L^1_2(\R^d)}\leq c(\|\widehat{C}\|_{L^1_4(\R^d)})\eps.
\eeq
Here we have also used $\|\eps \nu \chi_{B_{\eps^{-1/2}}}\widehat{C}\|_{L^1_2}\leq c\eps\|\widehat{C}\|_{L^1_3}$.

This allows us to estimate $\Bvechat_D-\Bvechat_{\text{app}}$. We have
$$
\begin{aligned}
\|\Bvechat_D-\Bvechat_{\text{app}}\|_{L^1(\R^d)}\leq & \eps^{1-d}\left\{\left\|\left(\widehat{D}\left(\frac{\cdot-k_0}{\eps}\right)-\widehat{C}\left(\frac{\cdot-k_0}{\eps}\right)\right)\vec{\eta}^{(j_0)}(k_0) \right\|_{L^1(B_{\eps^{r}}(k_0))}\right.\\
&+\left\|\widehat{D}\left(\frac{\cdot-k_0}{\eps}\right)\left(\vec{\eta}^{(j_0)}(\cdot)-\vec{\eta}^{(j_0)}(k_0)\right) \right\|_{L^1(B_{\eps^{r}}(k_0))}\\
&\left. +\left\|\widehat{C}\left(\frac{\cdot-k_0}{\eps}\right)\vec{\eta}^{(j_0)}(k_0)\right\|_{L^1(B^c_{\eps^{r}}(k_0))}  \right\}
\end{aligned}
$$
Next we use \eqref{E:D-C-est}, the Lipschitz continuity of $\vec{\eta}^{(j_0)}$, and the estimate $\|\widehat{C}\|_{L^1(B^c_{\eps^{r-1}})}\leq \eps^{(1-r)s_C}\|\widehat{C}\|_{L^1_{s_C}(\R^d)}$ for all $s_C\geq 0$. This produces at $r=1/2$ 
$$
\begin{aligned}
\|\Bvechat_D-\Bvechat_{\text{app}}\|_{L^1(\R^d)}\leq & c\left(\eps^2 + \eps^2\|\widehat{D}\|_{L^1_1(\R^d)}+ \eps^{1+s_C/2}\|\widehat{C}\|_{L^1_{s_C}(\R^d)} \right)\\
\leq & c(1+\|\widehat{D}\|_{L^1_1(\R^d)}+\|\widehat{C}\|_{L^1_4(\R^d)}) \eps^2,
\end{aligned}
$$
if $s_C=4$.

We can now summarize the error estimate
$$
\begin{aligned}
\|\Bvechat-\Bvechat_{\text{app}}\|_{L^1(\R^d)}&\leq \|\Bvechat_D-\Bvechat_{\text{app}}\|_{L^1(\R^d)} +\|\Bvechat-\Bvechat_D\|_{L^1(\R^d)} \\
&\leq\|\Bvechat_D-\Bvechat_{\text{app}}\|_{L^1(\R^d)} + \|\Bvechat_Q\|_{L^1(\R^d)}+ \|\Bvechat_R\|_{L^1(\R^d)}.
\end{aligned}
$$
The components $\Bvechat_Q$ and $\Bvechat_R$ are estimated in \eqref{E:BQ1-est}, \eqref{E:BQR-est}, and \eqref{E:BR-est}. Having now estimated $\|\widehat{R}\|_{L^1}$ in terms of $\|\widehat{C}\|_{L^1}$ and $\|\widehat{D}\|_{L^1}$ in terms of $\|\widehat{C}\|_{L^1_2}$, we get  for $r=1/2$ and $s_D=2$ 
$$\|\Bvechat_Q\|_{L^1(\R^d)}\leq c_1(\|\widehat{C}\|_{L^1_2(\R^d)}) \eps^{3}, \qquad \|\Bvechat_R\|_{L^1(\R^d)}\leq c_2(\|\widehat{C}\|_{L^1_2(\R^d)}) \eps^3,$$
where $c_1$ and $c_2$ depend polynomially on $\|\widehat{C}\|_{L^1_2(\R^d)}$.
Hence, the estimate in Theorem \ref{T:main} is proved.

\section{Numerical Example of Bifurcating Gap Solitons for $d=2$}
\label{S:num}
\begin{figure}[ht!]
\includegraphics[scale=0.6]{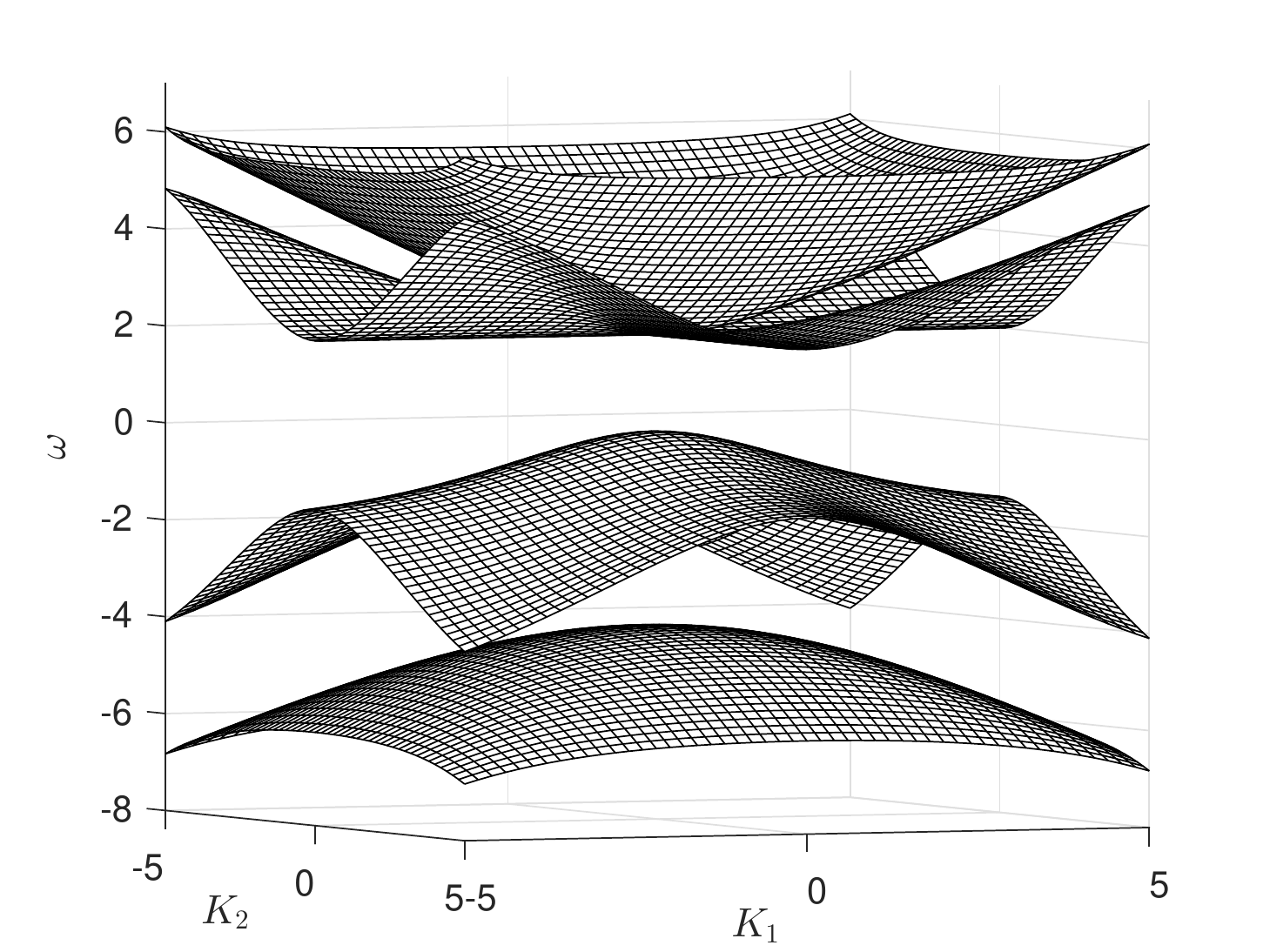}
\caption{Dispersion relation of \eqref{E:CME} for the example in Sec. \ref{S:num}.}
\label{F:disp_rel}
\end{figure}

In \cite{DW20} it is shown that assumption (A.1), i.e. the existence of a spectral gap is satisfied for $N=4$ in the symmetric case
\beq\label{E:lin-sym}
\begin{aligned}
&v_g^{(1)}=-v_g^{(2)}=:v,v_g^{(3)}=-v_g^{(4)}=:w, \\
&\kappa_{12}=\kappa_{34}=:\alpha_1,\\
&\kappa_{14}=\kappa_{32}=:\alpha_2,\\
&\kappa_{13}=\kappa_{42}=:\alpha_3,\\
&\kappa_{jj}=0, j=1,\dots,4 
\end{aligned}
\eeq
provided $|\alpha_1|^2>2(|\alpha_2|^2+|\alpha_3|^2)$.
In the following example we choose $v=(0,1)^T, w=(1,0)^T$, $\alpha_1=2,$ and $\alpha_2=\alpha_3=1$. The dispersion relation $\omega_j:\R^2\to\R, j=1,\dots,4$ of \eqref{E:CME} is plotted in Fig. \ref{F:disp_rel}. The gap appears even though the sufficient condition $|\alpha_1|^2>2(|\alpha_2|^2+|\alpha_3|^2)$ is not satisfied.
We see that the second eigenvalue $\lambda_2$ has an isolated maximum at $k=k_0:=0$. The corresponding frequency is $\omega_0:=\lambda_2(0)=0$. The eigenvector corresponding to $\lambda_2(0)$ is $\vec{\eta}_{j_0}(0)=\tfrac{1}{\sqrt{2}}(1,1,-1,-1)^T$. 

We use the following special case of the coefficients $\gamma_j^{(m,n,o)}$ in \eqref{E:B-eq}
\beq\label{E:nl-sym}
\begin{aligned}
1&=\gamma_{j}^{(j,j,j)}=\gamma_j^{(j,i,i)}=\gamma_j^{(i,i,j)}, \quad i,j=1,\dots,4,\\
&= \gamma_1^{(3,2,4)}= \gamma_1^{(4,2,3)}= \gamma_2^{(3,1,4)}= \gamma_2^{(4,1,3)}\\
&= \gamma_3^{(1,4,2)}= \gamma_3^{(2,4,1)}= \gamma_4^{(1,3,2)}= \gamma_4^{(2,3,1)},\\
&\gamma_j^{(m,n,o)}=0 \text{ otherwise.}
\end{aligned}
\eeq
Clearly, coefficients \eqref{E:lin-sym} and \eqref{E:nl-sym} allow symmetric solutions with $B_2=\overline{B_1}$ and $B_4=\overline{B_3}$. We do not make a direct use of this symmetry in our computations. We construct the approximation $\vec{B}_\text{app}$ of a solution of \eqref{E:B-eq} at $\omega = \omega_0+\eps^2\omega_1$ for six values of $\eps$: $0.2, 0.1, 0.05, 0.025, 0.0125$, and $0.00625$. 
The coefficients of the effective NLS \eqref{E:NLS} are $G_0=-0.25~ I_{2x2}$ and $\Gamma = 2.25$ and we choose $\omega_1=1$.
A real $C$ radially symmetric was chosen in this example. It was computed using the shooting method for the NLS in polar variables.

\begin{figure}[ht!]
\includegraphics[scale=0.6]{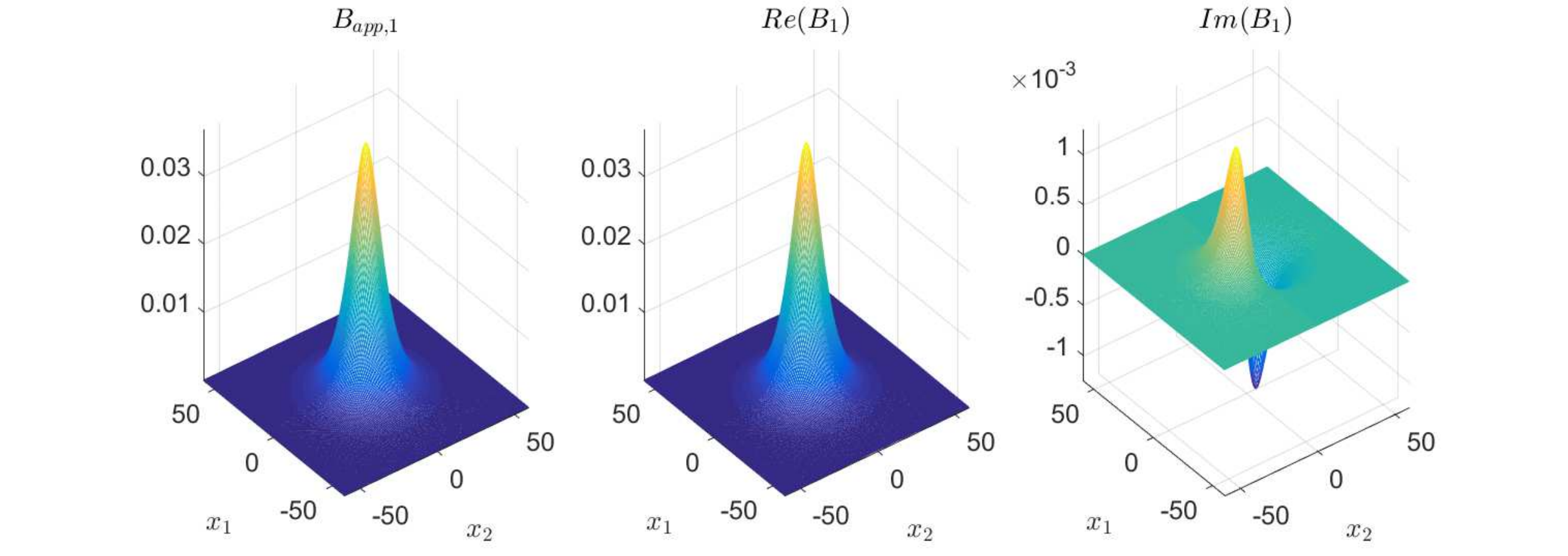}
\caption{Asymptotic approximation $B_{\text{app},1}$ and the numerical solution $B_1$ (real and imaginary part) at $\eps=0.05$.}
\label{F:profiles}
\end{figure}

Using the numerical Petviashvili iteration \cite{P76,AM05}, we also produce a numerical approximation of a solution $\vec{B}$ at $\omega = \omega_0+\eps^2$. The Petviashvili iteration is a fixed point iteration in Fourier variables with a stabilizing normalization factor. The initial guess of the iteration was chosen as $\vec{B}_\text{app}$. Note that although $\vec{B}_\text{app}$ can be real (if a real solution $C$ of the NLS is chosen), equation \eqref{E:B-eq} does not allow real solutions $\vec{B}$ due to the term $\ri \nabla \vec{B}$ and due to the realness of $\alpha_1,\alpha_2,$ and $\alpha_3$ and of $\gamma_j^{(m,n,o)}$. Nevertheless, if $\vec{B}_\text{app}$  is real, there must be a solution $\vec{B}$ with $\text{Im}(\vec{B})=O(\eps^2)$. Figure \ref{F:profiles} shows $\vec{B}_\text{app}$ and $\vec{B}$ for $\eps=0.05$. 

\begin{figure}[ht!]
\includegraphics[scale=0.6]{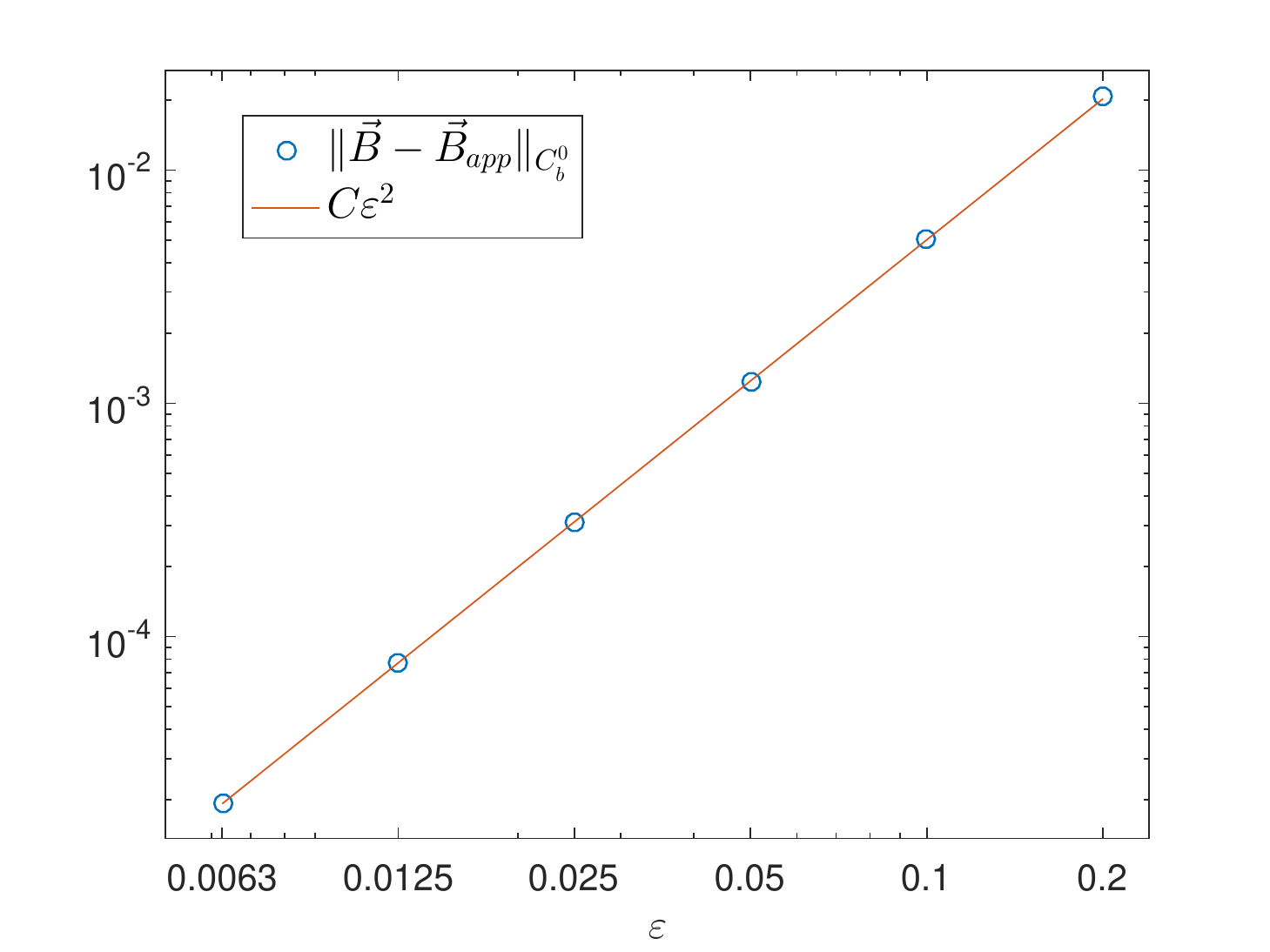}
\caption{Convergence of the asymptotic error in $\eps$. }
\label{F:eps_conv}
\end{figure}

The numerical parameters for the Petviashvili iteration were selected as follows: we compute on the domain $x\in [-3/\eps, 3/\eps]^2$ with the discretization given by 160x160 grid points, i.e. $dx_1=dx_2=3/(80\eps)$. Note that because $k_0=0$, the relatively coarse discretization for small values of $\eps$ does not matter (there are no oscillations to be resolved).
For each $\eps$ we evaluate the asymptotic error $E:=\|\vec{B}-\vec{B}_\text{app}\|_{C_0}$. Figure \ref{F:eps_conv} shows the convergence of the error in $\eps$. Clearly, $E(\eps)\sim c\eps^2$, which confirms the convergence rate proved in our theorem.

\section*{Acknowledgements}
This research is supported by the \emph{German Research Foundation}, DFG grant No. DO1467/3-1. 

\bibliographystyle{plain}
\bibliography{bib-CME-Rd}

\end{document}